\documentclass[12pt,twoside,a4paper,leqno]{article}
\usepackage[dvips,%
pdftitle={Finite order Vinogradov's C-spectral sequence},%
pdfauthor={R. Vitolo},%
pdfkeywords={Fibred manifold, jet space, variational bicomplex,%
variational sequence, spectral sequence}
pdfsubject={Finite order Vinogradov's C-spectral sequence},%
colorlinks,linkcolor={blue},citecolor={blue},urlcolor={red}%
]{hyperref}
\usepackage{amsmath}
\usepackage{amssymb}
\usepackage{amsthm}
\usepackage{diagrams}
\newcommand*{\email}[1]{\href{mailto:#1}{\begingroup \urlstyle{rm}\Url{#1}}}
\newcommand*{\eprint}[1]{\href{http://arXiv.org/abs/#1}%
{\begingroup \Url{arXiv:#1}}}
\newtheorem{thm}{THEOREM}
\newtheorem{cor}[thm]{COROLLARY}
\newtheorem{lem}[thm]{LEMMA}
\newtheorem{prop}[thm]{PROPOSITION}
\theoremstyle{definition}
\newtheorem{defn}[thm]{DEFINITION}
\newtheorem{REM}[thm]{REMARK}
\newcommand{\ie}{i.e$.$,}
\DeclareFontFamily{U}{UWCyr}{}
\DeclareFontShape{U}{UWCyr}{m}{n}{%
  <5> <6> <7> <8> <9>
  <10> <10.95> <12> <14.4> <17.28> <20.74> <24.88> wncyr10
  }{}
\DeclareMathAlphabet{\cyrm}{U}{UWCyr}{m}{n}
\DeclareSymbolFont{cyrm}{U}{UWCyr}{m}{n}
\DeclareSymbolFontAlphabet{\cyrm}{cyrm}
\newcommand*{\pd}[2]{\mathchoice{\frac{\partial #1}{\partial #2}}
  {\partial #1/\partial #2}{\partial #1/\partial #2} {\partial
  #1/\partial #2}}

\newcommand{\R}{\mathbb{R}}
\newcommand{\eval}[2][\right]{\relax
  \ifx#1\right\relax \left.\fi#2#1\rvert}

\let\abs=\envert

\def\QED{\hskip0.1em\hfill\null\
\null\nobreak\hfill\kern3pt\vbox{\hrule\hbox
   {\vrule\kern1pt\vbox{\kern1.7pt\hbox{$\scriptscriptstyle{QED}$}
    \kern0.2pt}\kern1pt\vrule}\hrule}}

\DeclareMathOperator{\byd}{\overset{\text{def}}{=}}
\DeclareMathOperator{\im}{im}

\DeclareMathOperator{\id}{id}

\DeclareMathOperator{\con}{\lrcorner}
\DeclareMathSymbol{\Evo}{\cyrm}{cyrm}{"03}
\DeclareMathOperator{\Vol}{Vol}
\newcommand{\cprime}{\/{\mathsurround=0pt$'$}}

\newcommand{\brh}{\boldsymbol{\rho}}
\newcommand{\bsi}{\boldsymbol{\sigma}}
\newcommand{\bta}{\boldsymbol{\tau}}

\newcommand{\hd}{\bar{d}}
\newcommand{\wid}{\widehat{d}}
\newcommand{\myskip}{\vspace*{8pt}}
\newcommand{\cC}{\mathcal{C}}
\newcommand{\cD}{\mathcal{D}}
\newcommand{\cE}{\mathcal{E}}
\newcommand{\cF}{\mathcal{F}}
\newcommand{\cV}{\mathcal{V}}
\newcommand{\cH}{\mathcal{H}}
\newcommand{\Diff}{\mathrm{Dif{}f}}
\DeclareMathOperator{\alt}{alt}

\DeclareMathOperator{\Hom}{Hom}
\newcommand{\hL}{\Bar{\Lambda}}
\newcommand{\CDiffalt}[2]{\cC\Diff^{\alt}_{(#1)\,#2}}
\newcommand{\CDiff}{\mathcal{C}\mathrm{Dif{}f}}
\title{Finite order formulation \\ of Vinogradov's $\cC$-spectral
    sequence} 
\author{Raffaele Vitolo
\\ {\footnotesize Dept. of Mathematics ``E. De Giorgi'', University of Lecce}
\\ {\footnotesize via per Arnesano, 73100 Lecce, Italy}
\\ {\footnotesize and Diffiety Institute, Russia}
\\ {\footnotesize email: \email{Raffaele.Vitolo@unile.it}
}}
\date{}
\begin{document}

\maketitle

\begin{abstract}
  The $\cC$-spectral sequence was introduced by Vinogradov in the late
  Seventies as a fundamental tool for the study of algebro-geometric
  properties of jet spaces and differential equations.  A spectral
  sequence arise from the contact filtration of the modules of forms on
  jet spaces of a fibring (or on a differential equation).  In order
  to avoid serious technical difficulties, the order of the jet space
  is not fixed, \ie{} computations are performed on spaces containing
  forms on jet spaces of any order.

  In this paper we show that there exists a formulation of Vinogradov's
  $\cC$-spectral sequence in the case of finite order jet spaces of a
  fibred manifold. We compute all cohomology groups of the finite
  order $\cC$-spectral sequence. We obtain a finite order variational
  sequence which is shown to be naturally isomorphic with Krupka's
  finite order variational sequence.


\end{abstract}
\section*{Introduction}

The framework of this paper is that of the algebro-geometric formalism
for differential equations.  This branch of mathematics started with
early works by Ehresmann and Spencer (mainly inspired by Lie), and was
carried out by several people.  One fundamental aspect is that a
natural setting for dealing with global properties of differential
equations (such as symmetries, integrability and calculus of
variations) is that of jet spaces (see~\cite{Kup80,MaMo83,Sau89} for
jets of fibrings and~\cite{KLV86,Vin84,Vin88} for jets of
submanifolds).  Indeed, any differential equation can be regarded as a
suitable submanifold of a jet space.

The \emph{$\cC$-spectral sequence} was introduced by Vinogradov in
the late Seventies~\cite{Vin77,Vin78,Vin84}.  It is a spectral
sequence (see~\cite{BoTu82,McC85} for a definition) arising from a
particular filtration of the De Rham complex on jet spaces (or on
differential equations).  Namely, there is a natural distribution on
any jet space, the \emph{Cartan distribution}, which consists of
tangent planes to the prolongation of any section of the jet space.
The space of forms annihilating the Cartan distribution, the
\emph{contact forms}, is an ideal of the space of all forms, and
yields a filtration of this space by means of its powers.  Contact
forms have deep meanings in several respects. For example, contact
forms yield zero contribution to action-like functionals. Indeed,
the $\cC$-spectral sequence yields the \emph{variational sequence}
as a by-product. The morphisms of the variational sequence are the
Euler--Lagrange morphism and other relevant maps from the calculus
of variations.

The above formulation has been carried out in the case of infinite
jets, in order to avoid serious technical difficulties due to the
computation of jet order. This paper provides a formulation of
Vinogradov's $\cC$-spectral sequence in the case of finite order jet
spaces of a fibred manifold. A finite order variational sequence then
arise from the formulation.  Of course, it was evident since the
earliest works that such a structure should exist. But the necessary
evaluations were not developed in general situations due to
severe technical difficulties.

As one could expect, the direct limit of the finite order
formulation yields the infinite order formulation by Vinogradov.
Moreover, a recent finite order variational sequence by Krupka has been
proposed~\cite{Kru90}; it is proved that our finite order formulation
recover Krupka's formulation, providing an equivalence between the
variational sequences of both cases.  As a by-product, it is possible
to represent any quotient space of the finite order variational sequence
by concrete sections of bundles. This problem is easily solved by
means of the intrinsic definition of adjoint operator \cite{Many99,Vin84}.

The equivalence between Vinogradov's and other infinite order
formulations~\cite{AnDu80,Bau82,OlSh78} either is evident or has
already been proved (see~\cite{Tul77,DeTu80} for a comparison between
Tulczyjew's and Vinogradov's formulation).  The advantage of the algebraic
techniques used in this paper is in the simplicity of their
generalisations.  For example, analogous result for jet spaces of
submanifolds of a given manifold could be obtained by a
straightforward generalisation.

Summarizing, this paper shows the possibility of computing the order
of objects involved at any step of the constructions even in the
infinite order formalism. It seems that working with infinite order
objects implies no loss of information because the order can always
be reconstructed. Moreover, working without a definite
order is easier than computing it every time. So, it seems that the
best strategy would be to compute it only when the problem being
investigated strictly requires it. Examples of such problems are
provided in \cite{Gri99,KrMu99}.

\myskip

We finish with some mathematical preliminaries.
\subsection*{Preliminaries}

In this paper, manifolds and maps between manifolds are $C^{\infty}$.
All morphisms of fibred manifolds (and hence bundles) will be
morphisms over the identity of the base manifold, unless otherwise
specified. All modules will be modules of sections of some vector bundles
(hence projective modules).

\myskip

Let $V$ be a vector space such that $\dim V = m$. Suppose that
$V=W_1\oplus W_2$, with $p_1:V\to W_1$ and $p_2:V\to W_2$ the related
projections. Then, we have the splitting
\begin{equation}
 \label{wed split}
   \wedge^{m}V=\bigoplus_{k+h=m}\wedge^{k}W_1\wedge\wedge^{h}W_2~,
\end{equation}
where $\wedge^{k}W_1\wedge\wedge^{h}W_2$ is the subspace of
$\wedge^{m}V$ generated by the wedge products of elements of
$\wedge^{k}W_1$ and $\wedge^{h}W_2$.

There exists a natural inclusion $\odot_{m}\,L(V,V) \subset
L(\wedge^{m}V, \wedge^{m}V)$. Then, the projections $p_{k,h}$ related
to the above splitting turn out to be the maps
\begin{displaymath}
  p_{k,h}=\binom{m}{k}\odot_{k}p_1\odot_{}\odot_{h}p_2
  :\wedge^{m}V\to\wedge^{k}W_1\wedge\wedge^{h}W_2.
\end{displaymath}

Let $V' \subset V$ be a vector subspace, and set $W'_{1} \byd
p_{1}(V')$, $W'_{2} \byd p_{2}(V')$.  Then we have
\begin{equation}\label{wed split and ssp}
V' \subset W'_{1} \oplus W'_{2},
\end{equation}
but the inclusion, in general, is not an equality.
\section{Jet spaces}

In this section we recall some facts on jet spaces.

Our framework is a fibred manifold
\begin{displaymath}
\pi : E \to M ,
\end{displaymath}
with $\dim M = n$ and $\dim E = n+m$.

We deal with the tangent bundle $TE \to E$, the tangent prolongation
$T\pi : TE \to T M$ and the vertical bundle $VE\byd \ker T\pi \to E$.

Moreover, for $0 \leq r$, we are concerned with the $r$-th jet space
$J^r\pi$; in particular, we set $J^0\pi \equiv E$. We recall the
natural fibrings
\begin{displaymath}
  \pi^r_s : J^r\pi \to J^s\pi,
  \qquad \pi^r : J^r\pi \to M,
\end{displaymath}
and the affine bundle $\pi^r_{r-1} : J^r\pi \to J^{r-1}\pi$, which is
associated with the vector bundle $\odot^r T^* M\otimes _{J^{r-1}\pi}
VE \to J^{r-1}\pi$ for $0 \leq s \leq r$. A detailed account of the
theory of jets can be found in~\cite{ALV91,Many99,MaMo83,Kup80,Sau89,Vin88}.

\myskip

Charts on $E$ adapted to the fibring are denoted by $(x^\lambda
,u^i)$.  Greek indices $\lambda,\mu,\dots$ run from $1$ to $n$ and
label base coordinates, Latin indices $i,j,\dots$ run from $1$ to $m$
and label fibre coordinates, unless otherwise specified.  We denote by
$(\pd{}{x^\lambda},\pd{}{u^i})$ and $(dx^\lambda,du^i)$,
respectively, the local bases of vector fields and $1$-forms on $E$
induced by an adapted chart.

We denote multi-indices by boldface Greek letters such as $\bsi =
(\sigma_1, \dots, \sigma_k)$, with $0 \leq \sigma_1, \dots,
\sigma_k\leq n$. We also set $|\bsi | \byd k$.

The charts induced on $J^r\pi$ are denoted by $(x^\lambda,u^i_{\bsi})$,
where $0 \leq |\bsi| \leq r$ and $u^i_{0} \byd u^i$.
The local vector fields and forms
of $J^r\pi$ induced by the fibre coordinates are denoted by
$(\pd{}{u_{\bsi}^i})$ and $(du^i_{\bsi})$, $0 \leq |\bsi| \leq r, 1
\leq i \leq m$, respectively.

A (local) section $s\colon M\to E$ can be prolonged to a section $j_rs\colon
M\to E$; if we set $u^i\circ s = s^i$, then we have the coordinate
expression
\begin{displaymath}
(j_rs)^i_{\bsi} = \pd{^{\abs{\bsi}}}{x^{\sigma_1}\dots \partial
  x^{\sigma_n}}s^i.
\end{displaymath}

A fundamental role is played by the \emph{contact maps} on jet spaces
(see~\cite{MaMo83}).  Namely, for $0 \leq r$, we consider the natural
inclusion over $J^{r}\pi$
\begin{displaymath}
  D_{r+1} :
  J^{r+1}\pi \to T^* M\otimes_{J^{r+1}\pi} TJ^r\pi.
\end{displaymath}
This inclusion comes from the fact that, if $s$ is a section of $\pi$,
then we can identify $j_{r+1}s$ with $Tj_rs$. Then, we have the
natural morphism
\begin{displaymath} \omega_{r+1} :
  J^{r+1}\pi \to T^*J^r\pi \underset{J^{r+1}\pi}{\otimes} VJ^r\pi,
\end{displaymath}
defined by $\omega_{r+1} = \id_{TJ^r\pi} - D_{r+1}$.  We have the
coordinate expressions
\begin{align*}
  D_r &= dx^\lambda\otimes D_\lambda = dx^\lambda\otimes
  \left(\pd{}{x^\lambda} + u^j_{\bsi\lambda}\pd{}{u^j_{\bsi}}\right),
  \\
  \omega_r &= \omega^j_{\bsi}\otimes\pd{}{u^j_{\bsi}} =
  (du^j_{\bsi}-u^j_{\bsi\lambda}dx^\lambda) \otimes\pd{}{u^j_{\bsi}},
\end{align*} for $0 \leq |\bsi| \leq r$.  We stress that
\begin{gather}
 \label{d e vart 1}
 D_r\con\omega_r=\omega_r\con D_r=0
 \\
 \label{d e vart 2}
 (\omega_r)^2=\omega_r\qquad\qquad(D_r)^2=D_r
\end{gather}
The (local) vector field ${D_r}_\lambda$ is said to be the
($\lambda$-th) \emph{total derivative operator}.

We can regard $D_{r+1}$ and $\omega^*_{r+1}$ as the injective fibred
morphism over $J^r\pi$
\begin{align*}
  &D_{r+1} : J^{r+1}\pi \underset{J^r\pi}{\times} T M \to J^{r+1}\pi
  \underset{J^r\pi}{\times} TJ^r\pi,
  \\
  &\omega_{r+1}^* : J^{r+1}\pi \underset{J^r\pi}{\times} V^*J^r\pi \to
  J^{r+1}\pi \underset{J^r\pi}{\times} T^*J^r\pi.
\end{align*}
We have the remarkable vector subbundles
\begin{align}
  \label{contact bundles}
  & C_{r+1,r} \byd \im D_{r+1} \subset
  J^{r+1}\pi\underset{J^r\pi}{\times}TJ^r\pi,
  \\
  & C^*_{r+1,r}\byd\im\omega_{r+1}^*\subset
  J^{r+1}\pi\underset{J^r\pi}{\times}T^*J^r\pi \subset T^*J^{r+1}\pi,
\end{align}
yielding the splitting~\cite{MaMo83}
\begin{equation}
 \label{jet connection}
 J^{r+1}\pi\underset{J^r\pi}{\times}T^*J^r\pi =\left(
   J^{r+1}\pi\underset{J^r\pi}{\times}T^* M\right) \oplus\im
 \omega_{r+1}^*.
\end{equation}

\myskip

Finally, we have a natural distribution $\cC_r$ on $J^r\pi$
generated by the tangent spaces to the prolongation $j_rs$ of any
section $s$, namely the \emph{Cartan distribution}
\cite{ALV91,Many99,MaMo83,Kup80,Sau89,Vin88}. It is generated by the
vector fields ${D_{r-1}}_\lambda$ and $\pd{}{u^i_{\bsi}}$, with
$\abs{\bsi}=r$. This distribution has not to be confused with
$C_{r,r-1}$, which is a subbundle of a different vector bundle (see
\eqref{contact bundles}), and is generated by $D_\lambda$.
\begin{REM}
  Both bundles $\cC_r$ and $C_{r,r-1}$ are part of chains of tangent
  projections whose inverse limit is the same, \ie{} Cartan
  distribution on infinite order jets, as it is immediate to show.
\end{REM}
\section{Vector fields and one--forms on jets}

Here we introduce distinguished modules over the ring of functions on
a jet space of a certain order. Namely, we denote by $\cF_{M}$ the
algebra $\cC^{\infty}(M)$, and by $\cF_r$ the algebra
$\cC^{\infty}(J^r\pi)$.

We denote by $\cD_{M}$ the module of vector fields on $M$, by $\cD_r$
the $\cF_r$-module (and $\R$-Lie Algebra) of vector fields on $J^r\pi$
and by $\cV_r\subset \cD_r$ the module of vertical vector fields.

It would be desirable to complement $\cV_r$ in $\cD_r$ with a natural
direct summand, like $\cC_r$. Unfortunately, this holds only in the
case of infinite order jet spaces, where Cartan distribution
provides the required summand.

We will encompass this problem by observing that another splitting
holds. Let us define a \emph{relative vector field} along
$\pi^{r+1}_r$ to be a map $X:J^{r+1}\pi \to TJ^r\pi$ such that
$\tau_{\pi^{r+1}_r} \circ X = \pi^{r+1}_r$. We denote by $\cD_{r+1,r}$
the $\cF_{r+1}$-module of relative vector fields along $\pi^{r+1}_r$.
In a similar way, we introduce the $\cF_{r+1}$-modules $\cC
\cD_{r+1,r}$ and $\cV_{r+1,r}$.

\begin{prop}
  We have the splitting
\begin{displaymath}
\cD_{r+1,r} = \cC  \cD_{r+1,r} \oplus \cV_{r+1,r}.
\end{displaymath}
The projections on the first and second factor are just contraction by
$\cD_{r+1}$ and $\omega_{r+1}$.
\end{prop}

We observe that $\cC \cD_{r+1,r}$ is locally generated over
$\cF_{r+1}$ by the sections ${D_{r+1}}_\lambda$.  Hence, any relative
vector field $X \in \cD_{r+1,r}$ can be split as
\begin{displaymath}
X = \Evo_X + \cC X,
\end{displaymath}
with the coordinate expression $\Evo_X = (X^i_{\bsi} -
u^i_{\bsi\lambda}X^\lambda)\pd{}{u^{\bsi}_i}$, and $\cC X = X^\lambda
D_\lambda$.  Pull-back yields the following inclusion
\begin{displaymath}
  {\pi^{r+1}_r}^*\cD_r \subset \cD_{r+1,r},
\end{displaymath}
so that the above splitting holds for any vector field $X \in \cD_r$.
In other words, $\Evo_X$ can be regarded as the \emph{vertical part}
and $\cC X$ as the \emph{horizontal part} of $X$ (see
\cite{ALV91,Many99,MaMo83,Kup80,Sau89,Vin88}).

\myskip

We consider the dual situation to the vector field case.

Let us set $\Lambda^{1}_r$ to be the $\cF_r$-module of $1$-forms on
$J^r\pi$. We introduce the submodule $\cH\Lambda^{1}_r \subset
\Lambda^{1}_r$ of forms with values in $T^* M$ (\emph{horizontal
  forms}) and the submodule $\cC^1\Lambda^{1}_r$ of forms
$\alpha\in\Lambda^{1}_r$ such that $(j_rs)^*\alpha=0$ for all sections
$s$ of $\pi$ (\emph{contact forms}). Note that the space of contact
forms is just the space of the annihilators of the Cartan distribution.

It would be desirable to complement $\cH\Lambda^1_r$ in
$\Lambda^{1}_r$ with a natural direct summand, like
$\cC^1\Lambda^{1}_r$. Unfortunately, this holds only in the case of
infinite order jet spaces, where the annihilator of the Cartan
distribution provides the required summand.

We will encompass this problem by observing that another splitting
holds.  Namely, we define $\Lambda^{1}_{r+1,r}$, $\cC^1\Lambda^{1}_{r+1,r}$
and $\cH\Lambda^{1}_{r+1,r}$ to be the $\cF_{r+1}$-modules of
$1$-forms on $J^{r+1}\pi$ with respective values in $T^*J^r\pi$,
$C_{r+1,r}^*$ and $T^* M$.

\begin{prop}\label{spl}
  We have the splitting
\begin{displaymath}
  \Lambda^{1}_{r+1,r} = \cC^1\Lambda^{1}_{r+1,r} \oplus
\cH\Lambda^{1}_{r+1,r}.
\end{displaymath}
The projections on the first and second factor are just contraction by
$\omega_{r+1}$ and $D_{r+1}$.
\end{prop}

If $\alpha\in\Lambda^{1}_{r+1,r}$ has the coordinate expression
$\alpha =\alpha_\lambda dx^\lambda +\alpha_i^{\bsi}du^i_{\bsi}$ ($0
\leq \bsi \leq r$), then
\begin{displaymath}
D_{r+1}(\alpha ) =(\alpha_\lambda+u^i_{\bsi\lambda}\alpha_i^{\bsi})\
dx^\lambda,
\qquad\quad \omega_{r+1}(\alpha ) =\alpha_i^{\bsi}\omega^i_{\bsi}.
\end{displaymath}
\section{Main splitting}

Here, for $k\leq 0$, we consider the standard $\cF_r$-module
$\Lambda^{k}_r$ of $k$-forms on $J^r\pi$ (which coincides with the
exterior power $\wedge^k\Lambda^1_r$).

For $p\leq k$ we introduce the ideal $\cC^p\Lambda^k_r$ of
$\Lambda^k_r$ generated by $p$-th exterior powers of
$\cC^1\Lambda^{1}_r$ (\emph{$p$-contact $k$-forms}). Of course,
$\cC^p\Lambda^{k}_r$ can be interpreted as the submodule of
$k$-forms which vanish when contracted with a $p$-vector with factors
in $\cC D_r(\pi)$.

We also introduce $\cH\Lambda^{k}_r$ (\emph{horizontal $k$-forms}) of
$\cH\Lambda^{1}_r$. Finally, we consider the obviously defined
$\cF_{r+1}$-modules $\Lambda^{k}_{r+1,r}$, $\cC^p\Lambda^k_{r+1,r}$ and
$\cH\Lambda^{k}_{r+1,r}$.

\begin{REM}
  We stress that pull-back via $\pi^{r+1}_r$ yields the inclusion
  $\Lambda^{k}_r \subset \Lambda^{k}_{r+1,r}$, and analogously for
  $\cC^p\Lambda^k_r$ and $\cH\Lambda^{k}_r$.
\end{REM}

If $\alpha \in \Lambda^{k}_{r+1,r}$, then we have
\begin{displaymath}
\alpha = \alpha
{_{i_1 \dots i_{h} }^{\bsi_1 \dots \bsi_{h}}}
{_{\lambda_{h+1} \dots \lambda_{k}}}
du^{i_1}_{\bsi_1}\wedge\dots\wedge du^{i_{h}}_{\bsi_{h}}\wedge
dx^{\lambda_{h+1}} \wedge\dots\wedge dx^{\lambda_{k}}.
\end{displaymath}
If $\beta\in\cH\Lambda^k_{r+1,r}$, then
\begin{displaymath}
\beta = \beta_{\lambda_{1} \dots \lambda_{k}} \, dx^{\lambda_{1}}
\wedge \dots \wedge dx^{\lambda_{k}}.
\end{displaymath}
If $\gamma\in\cC^k\Lambda^k_{r+1,r}$, then
\begin{displaymath}
\gamma = \gamma{_{i_{1}\dots i_{k}}^{\bsi_{1}\dots \bsi_{k}}}\,
   \omega^{i_{1}}_{\bsi_{1}} \wedge\dots\wedge
\omega^{i_{k}}_{\bsi_{k}},
    \qquad 0 \leq |\bsi_{1}|,\dots,|\bsi_{k}| \leq r.
\end{displaymath}
Here, the coordinate functions are sections of $\Lambda^{0}_{r+1}$,
and the indices' range is $0 \leq |\bsi_j| \leq r$, $0 \leq h \leq k$.
We remark that, in the coordinate expression of $\alpha$, the indices
$\lambda_j$ are suppressed if $h=k$, and the indices~${_{i_j}^{\bsi_j}}$
are suppressed if $h = 0$.

\myskip

In the rest of this section, we shall consider the effects of the
splitting of proposition~\ref{spl} on $\Lambda^k_r$. As one can
expect, $\Lambda^k_r$ does not split as a direct sum of exterior
products of $\cC^p\Lambda^k_{r}$ and $\cH\Lambda^{q}_{r}$, for
suitable $p$ and $q$. But we have the following result.

\begin{prop}\label{graded}
  The splitting of $\Lambda^{1}_{r+1,r}$ (proposition~\ref{spl})
  induces the splitting
\begin{displaymath}
\Lambda^{k}_{r+1,r}=\bigoplus_{l=0}^k
\cC^{k-l}\Lambda^{k-l}_{r+1,r}\wedge\cH\Lambda^{l}_{r+1}
\end{displaymath}
\end{prop}

We recall that, in the above splitting, direct summands with $l>n$
vanish.

\begin{defn}
  We define the above splitting to be the $\cC$-splitting.
\end{defn}

We set $p_{k-l,l}$ to be the projection of the $\cC$-splitting on the summand
$\cC^{k-l}\Lambda^{k-l}_{r+1,r}\wedge\cH\Lambda^{l}_{r+1}$. We set also
$H=p_{0,k}$ and $V=\id-H$. Due to results in the introduction, we have the
following theorem.

\begin{prop}
The explicit expression of the projections of the $\cC$-splitting is
\begin{displaymath}
  p_{k-l,l}=\binom{k}{k-l}\odot_{k-l}\omega_{r+1}\odot\odot_l D_{r+1}.
\end{displaymath}
\end{prop}

\begin{REM}\label{coord. expr. of h}
We have the coordinate expression
\begin{align*}
  & p_{k-l,l}(\alpha) = \sum u^{j_1}_{\bta_{1}\lambda_{1}} \dots
  u^{j_{s}}_{\bta_{s}\lambda_{s}} \alpha
  {_{i_1 \widehat{\dots} \, i_{k-l+s} \, j_{1} \dots j_{s}} ^{\bsi_1
      \widehat{\dots} \bsi_{k-l+s} \bta_{1} \dots \bta_{s}}}
  {_{\lambda_{s+1} \dots \lambda_{l}}}
  \\
  & \hphantom{h(\alpha) = \sum} \omega^{i_1}_{\bsi_1}
  \wedge\widehat{\dots}\wedge \omega^{i_{k-l+s}}_{\bsi_{k-l+s}} \wedge
  dx^{\lambda_{1}} \wedge\dots\wedge dx^{\lambda_{l}},
\end{align*}
where $0 \leq s \leq l$ and the sum is extended to all subsets
$\{ {^{j_1}_{\bta_{1}}} \dots {^{j_s}_{\bta_{s}}} \}$ of
$\{ {^{i_1}_{\bsi_1}} \dots {^{i_{k-l+s}}_{\bsi_{k-l+s}}} \}$,
where $\widehat{\dots}$ stands for suppressed indexes (and corresponding
contact forms) belonging to one of the above subsets.
\end{REM}

It turns out that, for $k \leq n$, we have the coordinate expression
\begin{displaymath}
H(\alpha) = u^{i_1}_{\bsi_{1}\lambda_{1}} \dots u^{i_{h}}_{\bsi_{h}\lambda_{h}}
\alpha
{_{i_1 \dots i_{h} }^{\bsi_1 \dots \bsi_{h}}}
{_{\lambda_{h+1} \dots \lambda_{k}}}
dx^{\lambda_{1}} \wedge\dots\wedge dx^{\lambda_{k}},
\end{displaymath}
with $0 \leq h \leq k$.

Now, we apply the conclusion of remark~\ref{wed split and ssp} of
introduction to the $\cF_r$-submodule $\Lambda^{k}_{r} \subset
\cH\Lambda^{k}_{r+1,r}$. To this aim, we want to find the image of
$\Lambda^{k}_{r}$ under the projections of the $\cC$-splitting. We need to
introduce further spaces.

We set $\cH_P\Lambda^{k}_{r+1}$ to be the $\cF_r$-module of horizontal
forms on $J^{r+1}\pi$ which are polynomials of degree $k$ with respect
to the affine structure of $J^{r+1}\pi \to J^r\pi$.

Moreover, we introduce the subspace
$\cC^k\Lambda^k_{r,r+1}\subset\Lambda^k_{r+1,r}$ of $k$-forms with values in
$C_{r+1}^*$ and coefficients in $\cF_r$.

Finally, we denote the restrictions of $H,V$ to $\Lambda^{k}_{r}$ by
$h,v$.

\begin{prop}\label{char. of h}
  Let $0 < k\leq n$, and denote
\begin{displaymath}
\hL^k_r \byd h (\Lambda^{k}_r).
\end{displaymath}

Then, we have the inclusion $\hL^k_r \subset
\cH_P\Lambda^{k}_{r+1}$.

Moreover, the $\cF_r$-module $\hL^{k}_r$ admits the following
characterization: $\alpha \in \cH_P\Lambda^{k}_{r+1}$ belongs to
$\hL^{k}_r$ if and only if there exists $\beta \in
\Lambda^{k}_{r}$ such that $(j_r s)^{*}\beta = (j_{r+1}s)^{*}\alpha$
for each section $s: M\to E$.
\end{prop}

\begin{proof}
  If $s : M \to E$ is a section, then the following identities
\begin{displaymath}
(j_{r}s)^*\beta = (j_{r+1}s)^*h(\beta),
\qquad\quad
(j_{r+1}s)^*v(\beta) = 0,
\end{displaymath}
yield
\begin{displaymath}
\alpha = h(\beta) \quad \Leftrightarrow \quad (j_{r}s)^*\beta =
(j_{r+1}s)^*\alpha
\end{displaymath}
for all $\alpha \in \cH_P\Lambda^{k}_{r+1}$ and $\beta \in
\Lambda^{k}_{r}$.
\end{proof}

\begin{REM}
  It comes from the above proposition that not any section of
  $\cH_P\Lambda^{k}_{r+1}$ is a section of $\hL^{k}_r$; indeed, a
  section of $\cH_P\Lambda^{k}_{r+1}$ in general contains \lq too many
  monomials\rq\ with respect to a section of $\hL^{k}_r$. This
  can be seen by means of the following example. Consider a one-form
  $\beta \in \Lambda^{1}_{0}$.  Then we have the coordinate
  expressions
\begin{displaymath}
\beta = \beta_{\lambda}d^\lambda + \beta_{i}d^i,
\qquad\quad
h(\beta) = (\beta_{\lambda} + u^i_\lambda\beta_{i}) d^\lambda.
\end{displaymath}
If $\alpha\in\cH_P\Lambda^{1}_{1}$, then we have the coordinate
expression $\alpha = (\alpha_{\lambda} +
u^j_{\mu}\alpha^{\mu}_{j}{_{\lambda}}) d^\lambda$.  It is evident
that, in general, there does not exists $\beta \in \Lambda^{1}_{0}$
such that $h(\beta) = \alpha$.
\end{REM}

\begin{cor}\label{Hor aff}
Let $\dim M = 1$. Then we have
\begin{displaymath}
\hL^{1}_r = \cH_P\Lambda^{1}_{r+1}.
\end{displaymath}
\end{cor}

\begin{lem}
The $\cF_{r+1}$-module morphisms $H,V$ restrict on $\Lambda^{k}_{r}$ to the
surjective $\cF_r$-module morphisms
\begin{displaymath}
h : \Lambda^{1}_{r} \to \hL^{1}_r,
\qquad
v : \Lambda^{1}_{r} \to \cC^1\Lambda^1_r.
\end{displaymath}
\end{lem}

\begin{proof}
The restriction of $H$ has already been studied.  As for the restriction of
$V$, it is easy to see by means of a partition of the unity that it is
surjective on $\cC^1\Lambda^1_r$.
\end{proof}

\begin{thm}\label{splitting fin.}
The $\cC$-splitting yields the inclusion
\begin{displaymath}
\Lambda^{k}_r \subset \bigoplus_{l=0}^{k}
\cC^{k-l}\Lambda^{k-l}_{r,r+1}\wedge\hL^{l}_r,
\end{displaymath}
and the splitting projections restrict to surjective maps.
\end{thm}

\begin{proof}
  In fact, for any $l \leq k$ the restriction of any projection of the
  $\cC$-splitting to $\Lambda^{k}_{r}$ is valued in the above spaces.
  Let us prove the surjectivity. Let $\Delta \in
  \cC^{k-l}\Lambda^{k-l}_{r,r+1}\wedge\hL^{l}_r$, where $0 \leq l
  \leq n$. We have the coordinate expression
\begin{align*}
  & \Delta = u^{j_1}_{\bta_{1}\lambda_{1}} \dots
  u^{j_{h}}_{\bta_{h}\lambda_{h}} \Delta
  {_{i_1 \, \dots \, i_{k-l} \, j_{1} \dots j_{h}} ^{\bsi_1 \dots
      \bsi_{k-l} \bta_{1} \dots \bta_{h}}} {_{\lambda_{h+1} \dots
      \lambda_{l}}}
  \\
  &\hphantom{ \Delta = u^{j_1}_{\bta_{1}\lambda_{1}} \dots
    u^{j_{h}}_{\bta_{h}\lambda_{h}} } \omega^{i_{1}}_{\bsi_{1}}
  \wedge\dots\wedge \omega^{i_{k-l}}_{\bsi_{k-l}} \wedge
  dx^{\lambda_{1}} \wedge\dots\wedge dx^{\lambda_{l}},
\end{align*}
where $0 \leq |\bsi_{i}|,|\bta_{i}| \leq r$ and $0 \leq h \leq l$.  If
$\{\psi_{i}\}$ is a partition of the unity subordinate to a coordinate
atlas of $ E$, then let
\begin{displaymath}
\alpha_{\Delta}{_i} \byd \psi_{i} \, \Delta
{_{t_1 \dots t_{r} }^{\brh_1 \dots \brh_{r}}}
{_{\lambda_{r+1} \dots \lambda_{k}}}
du^{t_{1}}_{\brh_{1}} \wedge\dots\wedge
du^{t_{r}}_{\bsi_{r}} \wedge
dx^{\lambda_{r+1}} \wedge\dots\wedge dx^{\lambda_{k}},
\end{displaymath}
where the set $\{ {^{t_1}_{\brh_{1}}} \dots
{^{t_r}_{\brh_{r}}} \}$ is a permutation of the set
$\{ {}^{i_1}_{\bsi_1} \dots {}
\text{\raisebox{-.5mm}{$^{i_{k-l}}_{\bsi_{k-l}}$}}
{}^{j_1}_{\bta_1} \dots {}^{j_l}_{\bta_l} \}$.  Then
$\sum_{i}\alpha_{\Delta}{_i} \in \Lambda^{k}_r$, and its projection on
$\cC^{k-l}\Lambda^{k-l}_{r,r+1}\wedge\hL^{l}_r$ is $\Delta$.

The proof is analogous for $k>n$.
\end{proof}

We remark that, in general, the above inclusion is a proper inclusion:
a sum of elements of the direct summands needs not to be an element
of $\Lambda^{k}_{r}$.

We have a final important consequence of the above results.

\begin{cor}\label{cC and ker h}
Let $p\leq k$. We have
\begin{displaymath}
\cC^1\Lambda^k_{r} = \ker h\, \quad \text{if} \quad 0\leq k \leq n,
\qquad
\cC^1\Lambda^k_{r} = \Lambda^k_{r} \quad \text{if} \quad k > n.
\end{displaymath}
\end{cor}

\begin{proof}
Let $\alpha \in \Lambda^{k}_{r}$. Then, for any section $s: M\to E$ we have
$(j_{r}s)^{*} \alpha = (j_{r+1}s)^{*} h(\alpha)$,
and $\alpha \in \ker h$ implies $\alpha \in \cC^p\Lambda^k_{r}$. Conversely,
suppose $\alpha \in \cC^p\Lambda^k_{r}$. Then we have
\begin{displaymath}
(j_{r+1}s)^{*} h(\alpha) =
h(\alpha)_{\lambda_{1}\dots\lambda_{k}} \circ j_{r+1}s \,\,\,
dx^\lambda_{1} \wedge\dots\wedge dx^\lambda_{k},
\end{displaymath}
hence $h(\alpha) = 0$.

The first assertion comes from the above identities and $\dim M =
n$.
\end{proof}
\section{Forms and differential operators}

The above construction could be reformulated in a purely algebraic
context (see, for example,~\cite{Kra97}). One important fact from this
theory is the \lq parallelism\rq{} between the language of forms and
the language of differential operators. This allows us to \lq
import\rq{} the theory of adjoint operators and Green's formula in our
setting.  To do this, we provide a natural isomorphism between the
module of contact forms and a space of differential operators.

\myskip

Let $P,Q$ be modules over an algebra $A$ over $\R$. We recall
(\cite{ALV91}) that a \emph{linear differential operator} of order $k$
is defined to be an $\mathbb{R}$-linear map $\Delta : P \to Q$ such
that
\[
[\delta_{a_0},[\dots ,[\delta_{a_k},\Delta]\dots]] =0
\]
for all $a_0, \dots ,a_k\in A$. Here, square brackets stand for
commutators and $\delta_{a_i}$ is the multiplication morphism. Of
course, linear differential operators of order zero are morphisms of
modules. The $A$-module of differential operators of order $k$ from
$P$ to $Q$ is denoted by $\Diff_k(P,Q)$. The $A$-module of
differential operators of any order from $P$ to $Q$ is denoted by
$\Diff(P,Q)$. This definition can be generalised to maps between the
product of the $A$-modules $P_1$,\dots $P_l$ and $Q$ which are
differential operators of order $k$ in each argument, \ie{}
\emph{multidifferential operators}. The corresponding space is denoted
by $\Diff_{k}(P_1,\dots,P_l;Q)$, or, if $P_1=\dots=P_l=P$, by
$\Diff_{(l)\, k}(P,Q)$. Accordingly, we define $\Diff_{(l)}(P,Q)$.

\myskip

When dealing with modules of sections over jets, it is convenient to
give a slightly more general definition of differential operator. In
fact, we have the natural inclusions $\cF_r\subset\cF_s$ for $r \leq
s$. So, if $P$ is a $\cF_r$-module and $Q$ is a $\cF_s$ module we can
introduce differential operators between $P$ and $Q$ in a natural way.
In particular, we consider operators whose expressions contain total
derivatives instead of standard ones. More precisely, we say a
differential operator $\Delta\colon P\to Q$ (of order $k$) to be
\emph{$\mathcal{C}$-dif\-fer\-en\-tial} if it can be restricted to the
manifolds of the form $j_r(M)$ and $j_s(M)$. In other words, $\Delta$ is a
$\mathcal{C}$-dif\-fer\-en\-tial operator if the equality
$j_rs(M)^*(\varphi)=0$, $\varphi\in P$, implies
$j_ss(M)^*(\Delta(\varphi))=0$ for any section $s\colon M\to E$.
In local coordinates, $\mathcal{C}$-dif\-fer\-en\-tial operators have
the form $(a^{\bsi}_{ij}D_{\bsi})$, where $a^{\bsi}_{ij}\in\cF_s$,
$D_{\bsi}=D_{\sigma_1}\circ\dots\circ D_{\sigma_k}$.

We denote the $\cF_s$-module of $\mathcal{C}$-dif\-fer\-en\-tial
operators of order $k$ from $P$ to $Q$ by $\CDiff_k(P,Q)$. We also
introduce the $\cF_s$ module of differential operators from $P$ to $Q$
of any order $\CDiff(P,Q)$. We can generalize the definition to
multi-$\mathcal{C}$-differential operators. In particular, we will be
interested to spaces of antisymmetric multi-$\mathcal{C}$-differential
operators, which we denote by $\CDiff^{\alt}_{(l)\,k}(P,Q)$.
Analogously, we introduce $\CDiff^{\alt}_{(l)}(P,Q)$.

\myskip

Next, we introduce a last important module of vector fields. Namely,
let us denote the $\cF_r$-module of relative vertical vector field
$\varphi\colon J^r\pi\to V\pi$ by $\varkappa_r$. Of course,
$\varkappa_0=\cV_0$. Then, any $\varphi\in\varkappa_r$ can be uniquely
prolonged to a relative vertical vector field $\Evo_\varphi\colon
J^{r+s}\pi\to V\pi_s$. It can be proved that $\varkappa$ contains all
non-trivial infinitesimal symmetries (even higher order ones) of the Cartan
distribution, see~\cite{Many99}, for example. If in coordinates
$\phi=\phi^i\pd{}{u^i}$, then
$\Evo_{\varphi}=D_{\bsi}\varphi^i\pd{}{u^i_{\bsi}}$. Such vector
fields are said to be \emph{evolutionary vector fields}.

\begin{prop}
We have the natural isomorphism
\begin{displaymath}
  \cC^p\Lambda^p_{r,r+1}\wedge\hL^{l}_r \to
  \CDiffalt{p}{r}(\varkappa_0,\hL^{l}_r),\qquad \alpha\mapsto\nabla_\alpha
\end{displaymath}
where $\nabla_\alpha(\varphi_1,\dots,\varphi_p)=
\Evo_{\varphi_1}\con(\dots\con(\Evo_{\varphi_p}\con\alpha)\dots )$.
\end{prop}

The above proposition can be proved by analogy with the infinite order
case (see~\cite{Many99}). Just recall that the isomorphism is realized
due to the fact that to any vertical tangent vector to $J^r\pi$ there
exists an evolutionary field passing through it.
\section{Spectral sequence}

The $\cC$-spectral sequence has been introduced by Vinogradov in the late
Seventies~\cite{Vin77,Vin78,Vin84}. It is a very powerful tool in the study of
differential equations and their symmetries and conservation laws.

Here, we present a new finite order approach to $\cC$-spectral
sequence on the jet space of order $r$ of a fibred manifold.
Such an approach has already been attempted in a
particular case~\cite{Duz83}. Indeed, the finite order formulation
presents some technical difficulties: our main tool is the splitting
of theorem~\ref{splitting fin.}, where the direct summands have a
rather complicated structure with respect to the infinite order analogue.
\subsection{Filtration}

The module $\Lambda^k_r$ is filtered by the submodules $\cC^p\Lambda^{k}_r$;
namely, we have the obvious \emph{finite} chain of inclusions
\begin{displaymath}
\Lambda^k_r \byd \cC^0\Lambda^{k}_r\supset \cC^1\Lambda^{k}_r \supset
\dots \supset
\cC^p\Lambda^{k}_r \supset \dots \supset \cC^I\Lambda^{k}_r \supset
\cC^{I+1}\Lambda^{k}_r = \{0\},
\end{displaymath}
where $I$ is the dimension of the contact distribution (see~\cite{Many99}).

\begin{defn}
We say the graded filtration
\begin{displaymath}
\{\cC^p\Lambda^{k}_r\}_{p \in \mathbb{N}}
\end{displaymath}
of $\Lambda^k_r$ to be the \emph{$\cC$-filtration} on the jet space of order
$r$.
\end{defn}

The $\cC$-filtration gives rise to a spectral sequence in a natural
way. The spectral sequences is a well-known tool in Algebraic Topology and
Homological Algebra (see, for example,~\cite{McC85}).

\begin{defn}
  We say the spectral sequence $(E^{p,q}_N, e_N)_{N,p,q \in
    \mathbb{N}}$ (with $p+q=k$) coming from the above filtration to be Vinogradov's
  \emph{$\cC$-spectral sequence of (finite) order $r$} on the fibred
  manifold $\pi$.
\end{defn}

The goal of next subsections is to describe all terms in the spectral
sequence that arise from the $\cC$-filtration.
\subsection{Spectral sequence: the term \protect$E_0$}

As a preliminary step for the study of $(E^{p,q}_0, e_0)$, we look for a
description of the spaces $\cC^p\Lambda^{k}_{r}$. To this aim, we introduce new
projections associated to the splitting of proposition~\ref{graded}.

Let $0 \leq q \leq n$; we denote by $H^p$ the projection
\begin{equation}\label{spectral sum}
\Lambda^{p+q}_{r+1,r} \to \bigoplus_{l=1}^{p}
\CDiffalt{p-l}{r}(\varkappa_0,\cH\Lambda^{q+l}_{r+1});
\end{equation}
we denote by $V^p$ the complementary projection, \ie{} $V^p=\id - H^p$.
Of course, $H^p = 0$ if $q=n$ and $H^1 = H$.  Also, we denote by $h^p$
and $v^p$ the corresponding restrictions to the subspace $\Lambda^{k}_{r}$.

\begin{REM}
By theorem~\ref{splitting fin.}, $h^p$ is not surjective on
$\oplus_{l=1}^{p}\cC^{p-l}_r\wedge\hL^{q+l}_r$ unless $p>1$ or
$q<n-1$ (in these cases there is only $1$ summand in~\eqref{spectral sum}).
\end{REM}

\begin{lem}\label{fin. ord. character. of C}
Let $p \geq 1$. Then, we have
\begin{align*}
&\cC^p\Lambda^{p+q}_{r} \simeq \ker h^p \quad \text{if} \quad q < n;
\\
&\cC^p\Lambda^{p+q}_{r} = \Lambda^{p+q}_{r} \quad \text{if} \quad q \geq n.
\end{align*}
\end{lem}

\begin{proof}
We recall that (corollary~\ref{cC and ker h}) the
theorem holds for $p=1$.
Then, we have the identities $\ker H^p = \im V^p$ and
$\im V^p = \langle (\im V)^p\rangle = \langle (\ker H)^p\rangle$, where
$\langle (\im V)^p \rangle$ denotes the ideal generated by $p$-th exterior
powers of elements of $\im V$ in $\Lambda^{1+q}_{r}$. So, by restriction to
$\Lambda^{k}_{r}$, we have $\ker h^p = \langle (\ker h)^p\rangle$. But, by
definition we have
$\cC^p\Lambda^{p+q}_{r} = \langle (\ker h)^p\rangle$, hence the result.
\end{proof}

Now, we compute $(E_{0},e_{0})$. We recall that
$E_{0}^{p,q} \equiv \cC^p\Lambda^{p+q}_{r} \big / \cC^{p+1}\Lambda^{p+q}_{r}$.
We denote also the differential $e_0$ (which is the quotient of $d$) by
$\hd\byd e_0$.

\begin{lem}
We have
\begin{align*}
&E_{0}^{p,0} = \ker h^p;
\\
&E_{0}^{p,q} \simeq \CDiffalt{p}{r}(\varkappa_0, \hL^{q}_r)
 \quad \text{if}\quad q \leq n;
\\
&E_{0}^{p,q} = \{0\} \quad \text{otherwise;}
\\
&\hd : E_{0}^{p,q}  \to
E_{0}^{p,q+1} : h^{p+1}(\alpha) \mapsto h^{p+2}(d\alpha).
\end{align*}
\end{lem}

\begin{proof}
The first and third assertions are trivial. As for the second one,
the inclusion is realized via the injective morphism
\begin{displaymath}
E_{0}^{p,q} \equiv \ker h^p \big / \ker h^{p+1} \to
\CDiffalt{p}{r}(\varkappa_0, \hL^{q}_r) \colon
[\alpha] \mapsto h^{p+1}(\alpha).
\end{displaymath}
The above morphism is also surjective: in fact, even if $h^{p+1}$ is
not surjective on its target space, it is surjective on each summand
of~\eqref{spectral sum}.

The differential $\hd$ can be read through the above morphism;
we obtain the last assertion.
\end{proof}

\begin{prop}
The bigraded complex $(E_{0},e_{0})$ is isomorphic to the sequence of
complexes
\newdiagramgrid{Vitolo}
{1.8,2.8,1.8,1.8,1.8}
{1,.5,1,.5,1,.5,1,.5,1,.3,1}
\begin{displaymath}
\begin{diagram}[grid=Vitolo]
0 & 0 & 0 & & 0
\\
\uTo_{\hd} & \uTo_{-\hd} & \uTo_{\hd} & & \uTo_{(-1)^I\hd}
\\
\hL^{n}_r &
\CDiffalt{1}{r}(\varkappa_0, \hL^{n}_r) &
\CDiffalt{2}{r}(\varkappa_0, \hL^{n}_r) &
\dots &
\CDiffalt{I}{r}(\varkappa_0, \hL^{n}_r)
\\
\uTo & \uTo & \uTo & & \uTo
\\
\dots & \dots & \dots & \dots & \dots
\\
\uTo_{\hd} & \uTo_{-\hd} & \uTo_{\hd} & & \uTo_{(-1)^I\hd}
\\
\hL^{1}_r & \CDiffalt{1}{r}(\varkappa_0, \hL^{1}_r) &
\CDiffalt{2}{r}(\varkappa_0, \hL^{1}_r) &
\dots &
\CDiffalt{I}{r}(\varkappa_0, \hL^{1}_r)
\\
\uTo_{\hd} & \uTo_{-\hd} & \uTo_{\hd} & & \uTo_{(-1)^I\hd}
\\
\Lambda^{0}_{r} & \CDiffalt{1}{r}(\varkappa_0, \cF_r) &
\CDiffalt{2}{r}(\varkappa_0, \cF_r) &
\dots &
\CDiffalt{I}{r}(\varkappa_0, \cF_r)
\\
\uTo & \uTo & \uTo & & \uTo
\\
0 & 0 & 0 & \dots & 0
\end{diagram}
\end{displaymath}

The sequence becomes trivial after the $I$-th column.

The minus signs are put in order to agree with an analogous convention on
infinite order variational bicomplexes.
\end{prop}

\begin{REM}\label{dv}
  The differential $\hd$ is different from the standard horizontal
  differential on jets~\cite{Many99,Sau89}, which we denote by $\wid$.
  In fact, $\wid$ is complemented to $d$ by the differential
  $\cC d$ according to $\pi^{r+1}_r{^*}d = \wid + \cC d$.
  This does not hold for $\hd$. For example, let $\alpha\in\hL^k_r$. Then
  \begin{align*}
    \hd\alpha &= h(d\alpha)
    \\
    & = h(\wid(h(\alpha)+v(\alpha)))+h(\cC d(h(\alpha)+v(\alpha)))
    \\
    & = h(\wid(h(\alpha)))+h(\wid(v(\alpha))),
  \end{align*}
  so that, in general, $\hd\neq\wid$ unless $v(\alpha)=0$.
\end{REM}
\subsection{Spectral sequence: the term \protect$E_1$}

In this section we describe the term $E_1$ of the spectral
sequence. We also show that the $\cC$-spectral sequence yields an exact
sequence of modules which is just the finite order version of the well-known
variational sequence.

We  recall that $E_{1} = H(E_{0})$, where the homology is taken
with respect to $\hd$. We start by determining the term $E_1^{p,n}$.

\begin{thm}\label{bicomplex}
We have the diagram
\newdiagramgrid{finitebicomplex}
{.8,.8,.8,1,1.4,1.5,1.5,1.5,1}
{1,.5,1,.5,1,.5,1}
\begin{displaymath}
\begin{diagram}[grid=finitebicomplex]
& & 0 & & & & 0 & &
\\
& & \uTo & & & & \uTo & &
\\
0 & \rTo & \hL^n_r \big/ \hd(\hL^{n-1}_r) &
\rTo^{e_1} & \dots & \rTo^{e_1} &
\CDiffalt{p}{r}(\varkappa_0,\hL^n_r)\big /\hd(E_{0}^{p,n-1}) &
\rTo^{e_1} & \dots
\\
& & \uTo & & & & \uTo & &
\\
& & \hL^n_r & & \dots & &
\CDiffalt{p}{r}(\varkappa_0, \hL^n_r) & & \dots
\\
& & \uTo_{\hd} & & & & \uTo_{(-1)^I\hd} & &
\\
& & \dots & & & & \dots & & \dots
\end{diagram}
\end{displaymath}
where the top row is a complex.
The bicomplex is trivial if $p>I$ and vertical arrows with values
into the quotients are trivial projections. We have the
identifications
\begin{align}
&E_{1}^{0,n} =
\hL^n_r \big / \hd(\hL^{n-1}_r),\label{functionals}
\\
&E_{1}^{p,n} =
\CDiffalt{p}{r}(\varkappa_0, \hL^n_r)\big /
\hd(\CDiffalt{p}{r}(\varkappa_0, \hL^{n-1}_r),\label{Euler and such}
\\
&e_{1}^{0,n} :
\hL^{n}_r \big / \hd(\hL^{n-1}_r) \to
(E_{0}^{1,n-1})\big /
\hd(\CDiffalt{1}{r}(\varkappa_0, \hL^{n-1}_r))\colon\notag
\\
&\hphantom{e_{1}^{p,1}\colon}
[h^{1}(\alpha)] \mapsto [h^{2}(d\alpha)],\notag
\\
&e_{1}^{p,n}\colon
\CDiffalt{p}{r}(\varkappa_0, \hL^{n}_r)\big /
\hd E_{0}^{p,n-1} \to
\CDiffalt{p+1}{r}(\varkappa_0, \hL^{n}_r)\big /
\hd E_{0}^{p+1,n-1}\colon\notag
\\
&\hphantom{e_{1}^{p,n}\colon}
[h^{p+1}(\alpha)] \mapsto [h^{p+2}(d\alpha)].\notag
\end{align}
\end{thm}

\begin{proof}
  The above identifications come directly from the definition of
  $E_{1}$.  As for the last statement, we have by definition (see,
  \emph{e. g.},~\cite{McC85})
\begin{displaymath}
e_{1}^{p,1} = \pi \circ \delta,
\end{displaymath}
where $\delta$ is the Bockstein operator induced by the exact sequence
and $\pi$ is the cohomology map induced by the corresponding map $\pi$
of the exact sequence.  So, suppose that
\begin{displaymath}
h^{p+1}(\alpha) \in E_{0}^{p,n} = \CDiffalt{p}{r}(\varkappa_0,\hL^{n}_r);
\end{displaymath}
we have $\alpha \in \Lambda^{p+n}_{r}$. Then,
\begin{displaymath}
\pi(d\alpha) = \hd(\pi(\alpha)) = 0,
\end{displaymath}
because $\hd$ raises the degree by $1$ on the horizontal factor, so,
$d\alpha \in \cC^{p+1}\Lambda^{p+1+n}_{r}$. Being $d(d\alpha) = 0$,
$d\alpha$ is closed in $\cC^{p+1}\Lambda^{p+1+n}_{r}$ under the
differential $d$, but $d\alpha$ is not exact in
$\cC^{p+1}\Lambda^{p+n}_{r}$, \ie{} there does not exist a form $\beta
\in \cC^{p+1}\Lambda^{p+n}_{r} = \ker h^{p+1}$ such that $d\beta =
\alpha$. Hence, $d\alpha$ determines a cohomology class $[d\alpha]$ in
$\cC^{p+1}\Lambda^{p+1+n}_{r}$ which is, by definition, the value of
$\delta([h^{p+1}(\alpha)])$. The map $\pi$ maps $d\alpha$ into
$h^{p+2}(d\alpha)$, hence the cohomology class $[d\alpha]$ is mapped
into $[h^{p+2}(d\alpha)]$ by $\pi$.
\end{proof}

\myskip

Now we determine $E_1^{p,q}$. We need some important preliminary
results.

We observe that the $\cC$-spectral sequence of order $r$ converges
to the de Rham cohomology of $E$. This is due to the fact that the
$\cC$-spectral sequence is a first quadrant spectral sequence. So,
according to the standard definition of convergence
\cite{BoTu82,McC85}, there exists $n_0\in\mathbb{N}$ such that
$E_{n_0}=E_s$ for $s > n_0$, and $E_{n_0}$ is isomorphic to the
quotient vector spaces of the filtration
\begin{displaymath}
  H^*(\Lambda^*)\supset iH^*(\cC^1\Lambda^*)\supset i^2H^*(\cC^2\Lambda^*)
  \supset \dots \supset i^IH^*(\cC^I\Lambda^*)\supset 0,
\end{displaymath}
of the de Rham cohomology of $E$ ($i$ is the natural inclusion).

\begin{lem}\label{le:exact contact}
  The sequence
  \newdiagramgrid{contact}
  {.8,.8,.8,.8,.8,.8,.8,.8,.8,.8,.8}
  {1}
  \begin{diagram}[grid=contact]
    0 & \rTo & \cC^p\Lambda_r^{p} & \rTo^d & \cC^p\Lambda_r^{p+1} & \rTo^d &
    \dots & \rTo^d & \cC^p\Lambda_r^{p+n-1} & \rTo^d & \dots
  \end{diagram}
  is exact up to the term $\cC^p\Lambda_r^{p+n-1}$.
\end{lem}
\begin{proof}
  We generalize arguments and computations from~\cite{Kru90}. The
  modules $\cC^p\Lambda_r^k$ are the spaces of global sections of the
  corresponding sheaves of contact forms. Such sheaves are soft
  sheaves because they are sheaves of modules over a sheaf of rings,
  $\cF_r$, which admit a partition of unity.  We want to prove that
  the corresponding sheaf sequence is exact up the term
  $\cC^p\Lambda_r^{p+n-1}$. In this case, such a sequence would be a soft
  resolution of $\cC^p\Lambda_r^{p}$, hence acyclic.

  First, we prove exactness at $\cC^p\Lambda_r^p$. We proceed by
  induction on $r$. Let
  $\alpha\in\cC^p\Lambda^p$ such that $d\alpha=0$. If we have the coordinate
  expression $\alpha=\alpha{_{i_{1}\dots i_{p}}^{\bsi_{1}\dots \bsi_{p}}} \,
   \omega^{i_{1}}_{\bsi_{1}} \wedge\dots\wedge
  \omega^{i_{p}}_{\bsi_{p}}$, where $\abs{\bsi_k}\leq r-1$, then
  \begin{equation}\label{eq:dalpha}
  \begin{split}
    d\alpha = d\alpha{_{i_{1}\dots i_{p}}^{\bsi_{1}\dots \bsi_{p}}} &\wedge
   \omega^{i_{1}}_{\bsi_{1}} \wedge\dots\wedge
  \omega^{i_{p}}_{\bsi_{p}}+
\\
   &\alpha{_{i_{1}\dots i_{p}}^{\bsi_{1}\dots \bsi_{p}}}
   \omega^{i_{1}}_{\bsi_{1}} \wedge\dots\wedge
   dx^\lambda\wedge\omega^{i_k}_{\bsi_k\lambda}\wedge\dots\wedge
  \omega^{i_{p}}_{\bsi_{p}}
  \end{split}
  \end{equation}
  where $1\leq k\leq p$. Hence
  $\alpha{_{i_{1}\dots i_{p}}^{\bsi_{1}\dots \bsi_{p}}}=0$ if
  $\abs{\bsi_k}= r-1$ for some $k$. The induction yields $\alpha=0$.

  Now, let $\alpha\in\cC^p\Lambda^h$, with $h>p$. We recall the
  \emph{contact homotopy operator}~\cite{Kru90}. Namely, let
  $(x^\lambda,y^i_{\bsi})$ be a fibred chart on $J^r\pi$. We define
  the map $H(t,x^\lambda,u^i_{\bsi})=(x^\lambda,tu^i_{\bsi}$), $0\leq
  t\leq 1$. We observe that
  $H^*\omega^i_{\bsi}=t\omega^i_{\bsi}+u^i_{\bsi}dt$. Moreover,
  $H^*\alpha=\alpha'\wedge dt +\alpha''$. So, we define the contact
  homotopy operator to be the map
  \begin{displaymath}
    A\alpha \byd \int_0^1\alpha'\wedge dt.
  \end{displaymath}
  It is easy to check that $\alpha=Ad\alpha+dA\alpha$. Now, let
  $d\alpha=0$. The proof is complete if we show that $A\alpha$ is a
  $p$-contact form. Unfortunately, the properties of $H^*$ imply that
  $A\alpha$ is in general a $(p-1)$-contact form. But $d\alpha=0$,
  hence we can prove that $\alpha = \beta+d\gamma$, where $\beta$ is
  $(p+1)$-contact and $\gamma$ is $p$-contact. We have the coordinate
  expression
  \begin{displaymath}
    \alpha = A{_{i_{1}\dots i_{p}}^{\bsi_{1}\dots \bsi_{p}}} \wedge
    \omega^{i_{1}}_{\bsi_{1}} \wedge\dots\wedge
    \omega^{i_{p}}_{\bsi_{p}} +
    d(B{_{i_{1}\dots i_{p}}^{\bsi_{1}\dots \bsi_{p}}} \wedge
    \omega^{i_{1}}_{\bsi_{1}} \wedge\dots\wedge
    \omega^{i_{p}}_{\bsi_{p}})
  \end{displaymath}
  where $A^{\dots}_{\dots},B^{\dots}_{\dots}\in\Lambda_r^{h-p}$ and
  $\abs{\bsi_k}\leq r-1$. Then the coordinate expression of
  $d\alpha$ is similar to \eqref{eq:dalpha}. Split
  $A^{\dots}_{\dots}=A^{\dots}_{\dots}{_h}+A^{\dots}_{\dots}{_c}$,
  where the first summand is an horizontal form and the second is a
  contact form. Then by induction on
  $r$ it is easy to see that $A^{\dots}_{\dots}{_h}=0$. Setting
  \begin{displaymath}
    \beta\byd  A{_{i_{1}\dots i_{p}}^{\bsi_{1}\dots \bsi_{p}}}_c \wedge
    \omega^{i_{1}}_{\bsi_{1}} \wedge\dots\wedge
    \omega^{i_{p}}_{\bsi_{p}}
  \end{displaymath}
  we have $\alpha=\beta+d\gamma$, $d\beta=0$, hence $\alpha=d(A\beta+\gamma)$.
\end{proof}

The above result allows us to compute the term $E_1^{p,q}$.

\begin{thm}\label{th:E1}
  We have
  \begin{enumerate}
  \item $E_1^{0,q}=H^q(E)$, for $q\neq n$;
  \item $E_1^{p,n}=H^p(E)$, for $p \geq 1$;
  \item $E_1^{p,q}=0$ for $q\neq n$ and $p\neq 0$.
  \end{enumerate}
\end{thm}
\begin{proof}
  The first result follows from the fact that $E_0^{0,q}$ is the
  quotient of the de Rham sequence with an exact sequence (see the
  above lemma), hence its
  cohomology is the de Rham cohomology of $J^r\pi$. This latter
  cohomology is equal to $H^*(E)$ because $J^r\pi$ has topologically
  trivial fibre over $E$.

  The third statement is a direct consequence of the above lemma

  The second statement comes from a straightforward computation and
  the convergence of the $\cC$-spectral sequence to the de Rham cohomology.
\end{proof}
\subsection{Spectral sequence: the variational sequence}

The results of theorem~\ref{th:E1} can be used to produce a new
sequence which is of fundamental importance, namely the variational sequence.
The spaces of the sequence are cohomology classes of $E_0$ and
$E_1$. In order to give an explicit expression to classes in $E_1$ a
key role will be played by the intrinsic definition of adjoint
operator~\cite{Many99}

\begin{thm}
We have the complex
\newdiagramgrid{varcomplex}
{.7,.7,2,2,3,1,1}
{1}
\begin{displaymath}
\begin{diagram}[grid=varcomplex]
\dots & \rTo^{\hd} & \hL^{n}_r & \rTo^{\Tilde{e}_1} &
\CDiffalt{1}{r}(\varkappa_0, \hL^{n}_r)\big /
\hd(\CDiff{1}{r}(\varkappa_0,\Lambda^{n-1}_r)) &
\rTo^{e_1} & \dots,
\end{diagram}
\end{displaymath}
where $\Tilde{e}_1$ is the map which make the following diagram
commute
\begin{displaymath}
  \begin{diagram}[tight]
     \hL^{n}_r & \rTo^{\Tilde{e}_1} & & &
     \CDiffalt{1}{r}(\varkappa_0, \hL^{n}_r)\big /
     \hd(\CDiffalt{1}{r}(\varkappa_0,\hL^{n-1}_r))
     \\
     & \rdTo & & \ruTo_{e_1} &
     \\
     & & \hL^n_r \big/ \hd(\hL^{n-1}_r) & &
  \end{diagram}
\end{displaymath}
The cohomology of the above complex turns out to be
naturally isomorphic to the de Rham cohomology of $E$.
\end{thm}

\begin{defn}
  We say the above complex to be the \emph{finite order $\cC$-variational
    sequence associated with the $\cC$-spectral sequence of order $r$
    on $\pi\colon E\to M$}.
\end{defn}

The word \lq variational\rq\ comes from the fact that we can identify
the objects of the space $\hL^{n}_r$ with \emph{$(r+1)$-st order
  Lagrangians}~\cite{Vit98,Vit99}. Moreover, next two spaces in the
sequence can be identified with a space of (finite order)
\emph{Euler--Lagrange morphisms} and a space of (finite order) \emph{Helmholtz
  morphism}, and the differential $e_1$ is the operator sending
Lagrangians into corresponding Euler--Lagrange morphism and
Euler--Lagrange type morphisms into Helmholtz morphisms.

\myskip

The $\cC$-variational sequence is defined through some quotient
spaces. Now, we prove that each equivalence class in these spaces can be
represented by a distinguished form.

To this aim, we observe that pull-back allows to take $\hd=\wid$.
More precisely, we can consider
$\alpha\in\hL^q_r\subset\cH\Lambda^q_{r+1}$, so that the
horizontalization on $(r+1)$-st order jets is the identity on
$\alpha$. In this way, the complementary map $v$ fulfills
$v(\alpha)=0$ on the $(r+1)$-st order jet. Hence
$\hd(\alpha)=\wid(\alpha)$.  Also, it is easily shown~\cite{Vit99}
that the $r$-th order $\cC$-variational sequence is embedded via
pull-back into the $(r+1)$-st order one. More generally, it can be
proved that the direct limit of the $r$-th order $\cC$-variational
sequence is just the standard infinite order $\cC$-variational
sequence (see also next subsection). Hence, we have the embedding
\begin{equation}\label{eq:embedding}
\begin{split}
  \CDiffalt{p}{r}(\varkappa_0, \hL^{n}_r)&\big /
     \hd(\CDiffalt{p}{r}(\varkappa_0,\hL^{n-1}_r)) \hookrightarrow
\\
  &\CDiffalt{p}{}(\varkappa, \hL^{n})\big /
     \hd(\CDiffalt{p}{}(\varkappa,\hL^{n-1})),
\end{split}
\end{equation}
where $\CDiffalt{p}{}$ is the space of operators of any order,
$\hL^{k}$ is the space of horizontal forms on the jet space of any
order, and $\varkappa$ is the space of relative vertical vector field
of any order.

Let $\cF$ be the space of functions on jet spaces of any order, and
set $\widehat{\varkappa}\byd\Hom_\cF(\varkappa,\hL^n)$.  We recall
that, if $\Delta\colon P\to Q$ is a $\cC$-differential operator,
then~\cite{Many99} there exists an operator $\Delta^*\colon
\widehat{Q}\to\widehat{P}$. It fulfills
\begin{equation}
  \label{eq:adjoint}
  \widehat{q}(\Delta(p))-(\Delta^*(\widehat{q}))(p) =
  \widehat{d}\omega_{p,\widehat{q}}(\Delta).
\end{equation}
In coordinates, if $\Delta=\Delta_{ij}^{\bsi}D_{\bsi}$, then
\begin{displaymath}
  \Delta^*=(-1)^{\abs{\bsi}}D_{\bsi}\circ\Delta_{ji}^{\bsi}.
\end{displaymath}

Now, the following well-known isomorphism holds (\cite{Vin77,Vin78,Vin84};
see also~\cite[p. 192]{Many99})
\begin{equation}\label{eq:repr}
  \CDiffalt{p}{}(\varkappa, \hL^{n})\big /
     \hd(\CDiffalt{p}{}(\varkappa,\hL^{n-1})) \simeq K_p(\varkappa),
\end{equation}
where  $K_p(\varkappa)\subset\CDiff^{\alt}_{(l-1)}
(\varkappa,\widehat{\varkappa})$ is the subspace of
operators $\nabla$ which are skew-adjoint in each argument, \ie{}
\begin{equation}\label{eq:repr2}
  (\nabla(\varphi_1,\dots,\varphi_{p-2}))^*=-\nabla(\varphi_1,\dots,
  \varphi_{p-2})
\end{equation}
for all $\varphi_1$, \dots, $\varphi_{l-2}\in \varkappa$.  Note that,
if $p=1$, then the isomorphism reads as the evaluation of the adjoint
of the given operator at the constant function $1$~\cite{Many99}.

The above considerations show that the equivalence class
$[\alpha]$ in the quotient space with contact degree $p$
is represented through the embedding~\eqref{eq:embedding} and the
isomorphism~\eqref{eq:repr} as the operator $\nabla_\alpha$ obtained
after skew-adjoining $\alpha$ in its first $(p-1)$-arguments and
adjoining it in its $p$-th argument.

Accordingly, the space of distinguished representatives of quotient
spaces of order $r$ is the subspace of $K_p(\varkappa)$ made by
$2^p(r+1)$-st order operators of the form of~\eqref{eq:repr2}.

Let us give a look in coordinates:
$\alpha\in\CDiffalt{1}{r}(\varkappa_0, \hL^{n}_r)$ has the expression
\begin{displaymath}
  \alpha = \alpha_{i_1\dots i_{p-1}j}^{\bsi_1\dots\bsi_{p-1}\bta}
  \omega^{i_1}_{\bsi_1}\wedge\dots\wedge
  \omega^{i_{p-1}}_{\bsi_{p-1}}\wedge\omega^j_{\bta}\wedge\Vol_M,
\end{displaymath}
where $\Vol_M\byd dx^1\wedge\dots\wedge dx^n$ is the local volume on
$M$. Hence, if $p=1$, then
\begin{displaymath}
  \nabla_\alpha= (-1)^{\abs{\bsi}}D_{\bsi}\alpha_i^{\bsi} \omega^i
  \wedge\Vol_M.
\end{displaymath}
This clearly shows that the first quotient space in the variational
sequence is the space of Euler--Lagrange type operators.
If $p=2$, then equations \eqref{eq:adjoint}, \eqref{eq:repr2} yield
\begin{equation}\label{eq:helm}
\begin{split}
  \nabla_{\alpha} =
  (-1)^{\abs{\bta}}\frac{1}{2}\left(
  D_{\bta}\alpha_{ij}^{\bsi\bta}
  -\sum_{\abs{\brh}=0}^{s-\abs{\bsi}}(-1)^{\abs{(\bsi,\brh)}}
  \binom{\abs{(\bsi,\brh)}}{\abs{\brh}}
  D_{\brh}D_{\bta}\alpha_{ji}^{(\bsi,\brh)\bta}
  \right)
\\
  \omega^i_{\bsi}\wedge\omega^j\wedge\Vol_M,
\end{split}
\end{equation}
where $(\bsi,\brh)$ denotes the union of the multiindexes $\bsi$ and
$\brh$, $s$ is the jet order of $D_{\bta}\alpha_{ij}^{\bsi\bta}$,
and the factor $1/2$ comes from skew-symmetrization. Note that we
also used the Leibniz rule for total derivatives~\cite{Sau89} to derive
the expression of the $\cC$-differential operator
$(-1)^{\bsi}D_{\bsi}\circ\beta_{ji}^{\bsi}$, with
$\beta_{ji}^{\bsi}= D_{\bta}\alpha_{ij}^{\bsi\bta}$.

This clearly shows that the second quotient space in the variational
sequence is the space of Helmholtz type operators. We recall that
the Helmholtz operator of an Euler--Lagrange type operator $\eta$ is
just $e_1(\eta)$. If $e_1(\eta)=0$ then the local exactness of the
variational sequence tells us that $\eta$ is (locally) the
Euler--Lagrange operator of a Lagrangian.

Through the above expressions it is possible to derive a
representation formula for any $p$.

\begin{REM}\label{rem:e1}
  The expression of $e_1$ between quotient spaces can be derived by
  observing that $e_1$ is the quotient of the contact (or vertical)
  differential $\cC d$~\cite{Vit98,Vit99}. Namely, $e_1([\alpha])= [\cC
  d(\alpha)]$.
\end{REM}

\begin{REM}
  We could consider the \lq complementary\rq{} problem to the
  representative's one. More precisely, given
  $\alpha\in\CDiffalt{p}{r}(\varkappa_0, \hL^{n}_r)$ we can look for a
  section $q$ fulfilling $\alpha = \nabla_\alpha + \hd q$. Such a
  section always exists due to the vanishing of the cohomology of
  $\hd$ on the space where $q$ lives.

  This problem, for $p=1$, is the search for a Poincar\'e--Cartan form
  (see~\cite{Kol83,Vit98}, for example). Alonso Blanco~\cite{AB01}
  is able to determine one section $q$ through a connection for any $p$.
\end{REM}
\subsection{Comparisons}
\label{sec:comparisons}

We can perform two different kinds of comparisons.

\myskip

\textbf{Comparison with the standard infinite order $\cC$-spectral
  sequence.} The analysis has been done in~\cite{Vit99}. Here we
summarize the main steps.

It is clear that pull-back provides an inclusion of the $r$-th order
$\cC$-spectral sequence into the $(r+1)$-st order $\cC$-spectral
sequence. This yields a sequence of spectral sequences whose direct
limit is the infinite order $\cC$-spectral sequence~\cite{Vit99}.
By the way, it can be easily proved (see remark~\ref{dv}) that
the direct limit of $\hd$ and $\wid$ is the same.

One of the main differences between the finite order approach and the infinite
order approach is that in the infinite order case the diagram
\eqref{bicomplex} is a bicomplex, where horizontal arrows are provided
by $\cC d$ (see remark~\ref{dv}). In particular, it is easily proved
that the differential $e_1$ turns out to be the equivalence class of
the differential $\cC d$. In the finite order case this differential
does not provide any additional complex because it raises the order of
jet by one.

Another important point is that in the finite order approach there are
only finitely many non-zero $E_N^{p,q}$, because de Rham complex on
finite order jets \lq stops\rq .

In the finite order approach we recover some well-known results
of the infinite order case but proofs are slightly different. This can
give some new insight in the theory.

\myskip

\textbf{Comparison with Krupka's finite order variational sequence.}
Krupka's sequence~\cite{Kru90} is obtained by quotienting the de Rham
sequence on $J^r\pi$ with a natural subsequence. The first $n$ terms
of this subsequence are given in lemma~\ref{le:exact contact}, next
ones are the spaces of sections which are \emph{locally} in the space
$\ker h^p+d\ker h^{p-1}$. In~\cite{Vit98,Vit99} it is proved that
Krupka's sequence is isomorphic to the sequence
\begin{gather*}
\begin{diagram}
\dots & \rTo^{\cE_{n-1}} &
\hL^n_r & \rTo^{\cE_n} &
\CDiffalt{1}{r}(\varkappa_0, \hL^n_r)\big / h^{2}(d\ker h^1)^s &
\rTo^{\cE_{n+1}} &
\end{diagram}
\\
\begin{diagram}
\dots & \rTo^{\cE_{n+i-1}} &
\CDiffalt{i}{r}(\varkappa_0, \hL^n_r) \big / h^{i+1}(d\ker h^i)^s &
\rTo^{\cE_{n+i}} &
\dots
\end{diagram}
\end{gather*}
where $\cE_{k}([h^p(\alpha)]) = [h^{p+1}(d\alpha)]$, and
$h^{i+1}(d\ker h^i)^s$ stands for the space of sections which are
\emph{locally} in the space $h^{i+1}(d\ker h^i)$. Moreover, in~\cite{Vit99} it
is proved that the above sequence and the finite order
$\cC$-variational sequence are isomorphic up to the $(n+1)$-st term.

We can easily complete the proof and show the equivalence of both finite
order sequences in all their terms. In fact,
\begin{displaymath}
  h^{i+1}(d\ker h^i)^s=h^{i+1}(d\ker h^i),
\end{displaymath}
because the cohomology of $e_0=\hd$ vanish (see theorems
\ref{bicomplex} and~\ref{th:E1}).

In particular, the second quotient space of the $\cC$-variational
sequence is naturally isomorphic to the $(n+2)$-th space of Krupka's
variational sequence. This means that we provided in equation
\eqref{eq:helm} a way to represent all elements of degree $(n+2)$ of Krupka's
variational sequence, and the way to a possible generalization to $(n+p)$.

\myskip

\textbf{Comparison with Anderson and Duchamp's approach.} There is
only one more approach to finite order variational
sequences~\cite{AnDu80}, but the sequence provided in that paper
stops after the $(n+1)$-st term. Moreover, this sequence has been
derived directly from coordinate computations. Also, from the
coordinate expression it is easy to see that the sequence coincides
with the $\cC$-variational sequence up to the $n$-th term. The last
map in Anderson and Duchamp's paper has the same value of the
corresponding map in the $\cC$-variational sequence but is defined
between slightly different spaces.

To the author's knowledge the $\cC$-variational sequence, Krupka's
variational sequence and Anderson and Duchamp's sequence are the only
finite order formulation of variational sequences.
\section{Conclusions}

We derived a new finite order formulation of Vinogradov's $\cC$-spectral
sequence. We recovered most results of the infinite order formulations
in the finite order case. We have shown that the associated
$\cC$-variational sequence is isomorphic to Krupka's one. But the
techniques employed throughout this paper can be easily generalized,
e.g., to jets of submanifolds and differential equations, along the
lines of similar developments in the infinite order case. Moreover, we
have briefly seen how the infinite order scheme can be reached from
the finite order one through a direct limit process.

It seems to be clear that working on finite order jets produces almost
the same amount of information as working on infinite jets, from a \lq
structural\rq\ viewpoint. But there is a cost in carrying the order of
jets in all computations.

The author's opinion is that it is possible to skip problems related
to the order when the target of the research does not involve it.
Actually, in this paper we have shown that there exists the
possibility of computing orders also in the infinite order
$\cC$-spectral sequence. It is not necessary to use such a
possibility from the very beginning of any investigation: in fact,
the order can be computed at any step. Of course, there are problems
where the order is fundamental. Main examples are the problem of the
minimal order variational potential for a variationally trivial
Lagrangian and the problem of the minimal order Lagrangian for a
given variationally trivial Euler--Lagrange operator
\cite{Gri99,KrMu99}.

\subsection*{Acknowledgements}
  I would like to thank J. S.  Krasil{\cprime}shchik, M.  Modugno, and
  A. M. Vinogradov for discussions and comments. Special thanks are
  also due to A. M. Verbovetsky for his enlightening explanations about
  secondary calculus. This paper has been written during and after
  the stimulating courses of the \emph{Diffiety School} (Forino,
  Italy, July 1998) and the \emph{School in Homological Methods in
    Equations of Mathematical Physics} (Levo\v ca, Slovakia, August
  1998).

  This paper has been partially supported by GNFM of INdAM, MURST,
  Universities of Florence and Lecce.

\end{document}